\documentclass[12pt,twoside]{amsart}

\begin{document}

\centerline{}

\centerline {\Large{\bf Intuitionistic fuzzy $\Psi$-$\Phi$-contractive mappings }}
\centerline{\Large{\bf and fixed point theorems in non-Archimedean }}
\centerline{\Large{\bf intuitionistic fuzzy metric spaces}}
\centerline{}
%% My definition
\newcommand{\mvec}[1]{\mbox{\bfseries\itshape #1}}

\centerline{\bf {B. Dinda$^1$, T.K. Samanta$^2$ and Iqbal H. Jebril$^3$}}

\centerline{}

\centerline{$^1$Department of Mathematics,}
\centerline{Mahishamuri Ramkrishna
Vidyapith,}
\centerline{Howrah-711401, West Bengal, India. }
\centerline{E-mail: bvsdinda@gmail.com}
\centerline{}
\centerline{$^2$Department of Mathematics,}
\centerline{Uluberia College, Howrah, India.}
\centerline{E-mail: mumpu$_{-}$tapas5@yahoo.co.in}
\centerline{}
\centerline{$^3$Department of Mathematics,}
\centerline{King Faisal University,}
\centerline{Kingdom of Saudi Arabia}
\centerline{E-mail:ijopcm@yahoo.com}
\centerline{}

\newtheorem{Theorem}{\quad Theorem}[section]

\newtheorem{definition}[Theorem]{\quad Definition}

\newtheorem{theorem}[Theorem]{\quad Theorem}

\newtheorem{remark}[Theorem]{\quad Remark}

\newtheorem{corollary}[Theorem]{\quad Corollary}

\newtheorem{note}[Theorem]{\quad Note}

\newtheorem{lemma}[Theorem]{\quad Lemma}

\newtheorem{example}[Theorem]{\quad Example}

\centerline{}
\centerline{\bf Abstract}
\textbf{\emph{In this paper intuitionistic fuzzy $\psi$-$\phi$ contractive mappings are introduced. Intuitionistic fuzzy Banach contraction theorem for M-complete non-Archimedean intuitionistic fuzzy metric spaces and intuitionistic fuzzy Elelstein contraction theorem for non-Archimedean intuitionistic fuzzy metric spaces by intuitionistic fuzzy $\psi$-$\phi$ contractive mappings are proved.}}\\


%=============================
\section{\bf Introduction}
%=============================
Theory of intuitionistic fuzzy set as a generalization of fuzzy set \cite{zadeh} was introduce by Atansov \cite{Atanassov}. Grabiec \cite{grabiec} initiated the study of fixed point theory in fuzzy metric spaces. George and Veeramani \cite{veeramani} have pointed out that the definition of Cauchy sequence for fuzzy metric spaces given by Grabiec \cite{grabiec} is weaker and they gave one stronger definition of Cauchy sequence and termed as M-Cauchy sequence. The definition of Cauchy sequence given by Grabiec \cite{grabiec} has been termed as G-Cauchy sequence. With the help of fuzzy $\psi$-contractive mappings defined by Dorel Mihet\cite{mihet}, we introduce intuitionistic fuzzy $\psi$-$\phi$ contractive mappings. Our definition of intuitionistic fuzzy $\psi$-$\phi$ contractive mapping is more general than the definitions of intuitionistic fuzzy contractive mapping given by Abdul Mohamad \cite{abdul} and by this contraction we prove an intuitionistic fuzzy Banach contraction theorem for M-complete non-Archimedean intuitionistic fuzzy metric spaces. We also prove an intuitionistic fuzzy Elelstein contraction theorem for non-Archimedean intuitionistic fuzzy metric spaces.\\

%==============================
\section{\bf Preliminaries}
%==============================
We quote some definitions and statements of a few theorems
which will be needed in the sequel.

\begin{definition}
\cite{Schweizer} A binary operation \, $\ast \; : \; [\,0 \; , \;
1\,] \; \times \; [\,0 \; , \; 1\,] \;\, \rightarrow \;\, [\,0
\; , \; 1\,]$ \, is continuous \, $t$ - norm if \,$\ast$\, satisfies
the
following conditions \, $:$ \\
$(\,i\,)$ \hspace{0.5cm} $\ast$ \, is commutative and associative ,
\\ $(\,ii\,)$ \hspace{0.4cm} $\ast$ \, is continuous , \\
$(\,iii\,)$ \hspace{0.2cm} $a \;\ast\;1 \;\,=\;\, a \hspace{1.2cm}
\forall \;\; a \;\; \varepsilon \;\; [\,0 \;,\; 1\,]$ , \\
$(\,iv\,)$ \hspace{0.2cm} $a \;\ast\; b \;\, \leq \;\, c \;\ast\; d$
\, whenever \, $a \;\leq\; c$  ,  $b \;\leq\; d$  and  $a \,
, \, b \, , \, c \, , \, d \;\, \varepsilon \;\;[\,0 \;,\; 1\,]$.
\end{definition}
A few examples of continuous t-norm are $\,a\,\ast\,b\,=\,ab,\;\,a\,\ast\,b\,=\,\min\{a,b\},\;\,a\,\ast\,b\,=\,\max\{a+b-1,0\}$.

\begin{definition}
\cite{Schweizer}. A binary operation \, $\diamond \; : \; [\,0 \; ,
\; 1\,] \; \times \; [\,0 \; , \; 1\,] \;\, \rightarrow \;\,
[\,0 \; , \; 1\,]$ \, is continuous \, $t$-conorm if
\,$\diamond$\, satisfies the
following conditions \, $:$ \\
$(\,i\,)\;\;$ \hspace{0.1cm} $\diamond$ \, is commutative and
associative ,
\\ $(\,ii\,)\;$ \hspace{0.1cm} $\diamond$ \, is continuous , \\
$(\,iii\,)$ \hspace{0.1cm} $a \;\diamond\;0 \;\,=\;\, a
\hspace{1.2cm}
\forall \;\; a \;\; \in\;\; [\,0 \;,\; 1\,]$ , \\
$(\,iv\,)$ \hspace{0.1cm} $a \;\diamond\; b \;\, \leq \;\, c
\;\diamond\; d$ \, whenever \, $a \;\leq\; c$  ,  $b \;\leq\; d$
 and  $a \, , \, b \, , \, c \, , \, d \;\; \in\;\;[\,0
\;,\; 1\,].$
\end{definition}
A few examples of continuous t-conorm are $\,a\,\diamond\,b\,=\,a+b-ab,\;\,a\,\diamond\,b\,=\,\max\{a,b\},\;\,a\,\diamond\,b\,=\,\min\{a+b,1\}$.

\begin{definition}\cite{park}
A 5-tuple $(X,\mu,\nu,\ast,\diamond)$ is said to be an intuitionistic fuzzy metric space if $X$ is an arbitrary set, $\ast$ is a continuous t-norm, $\diamond$ is a continuous t-conorm, $\mu$ and $\nu$ are fuzzy sets on $X^2 \times (0,\infty)\,$ and $\,\mu$ denotes the degree of nearness, $\nu$ denotes the degree of non-nearness between $x$ and $y$ relative to $t$ satisfying the following conditions: for all $x,y,z\in\,X,\,s,t>0,\\\\$
$(\,i\,)$ \hspace{0.10cm}  $\mu(x,y,t) \,+\, \nu(x,y,t) \,\leq\, 1 $\\
$(\,ii\,)$ \hspace{0.10cm}$\mu(x,y,t) \,>\, 0 \, ;$ \\
$(\,iii\,)$ $\mu(x,y,t)\,=\, 1\,$ if and only if \, $x \,=\, y \,$ ; \\
$(\,iv\,)$\hspace{0.05cm} $\mu(x,y,t)=\mu(y,x,t)$ ;\\
$(\,v\,)$ \hspace{0.10cm} $\mu(x,z,t+s)\,\geq\,\mu(x,y,t)\,\ast\,\mu(y,z,s)$ ;\\
$(\,vi\,)$ \hspace{0.05cm}$ \mu(x,y,\cdot)\,:\,(0,\infty)\,\rightarrow\,(0,1]$ is continuous;\\
$(\,vii\,)$ \hspace{0.10cm}$\nu(x,y,t) \,>\, 0 \, ;$\\
$(\,viii\,)$ $\nu(x,y,t)\,=\, 0\,$ if and only if \, $x \,=\, y \,$ ; \\
$(\,ix\,)$\hspace{0.05cm} $\nu(x,y,t)=\nu(y,x,t)$ ;\\
$(\,x\,)$ \hspace{0.15cm} $\nu(x,z,t+s)\,\leq\,\nu(x,y,t)\,\diamond\,\nu(y,z,s)$ ;\\
$(\,xi\,)$ \hspace{0.04cm}$ \nu(x,y,\cdot)\,:\,(0,\infty)\,\rightarrow\,(0,1]$ is continuous.
\end{definition}

\begin{remark}
If in the above definition the triangular inequalities $(v)$ and $(x)$ are replaced by\\
$\mu(\;x,\,z\,,\,\max\{t,s\}\;)\,\geq\,\mu(x,y,t)\,\ast\,\nu(y,z,s)\;\,$ and\\
$\nu(\;x,\,z\,,\,\max\{t,s\}\;)\,\leq\,\nu(x,y,t)\,\diamond\,\nu(y,z,s).$\\
Or, equivalently,\\
$\mu(x,z,t)\,\geq\,\mu(x,y,t)\,\ast\,\mu(y,z,t)\;\,$ and\\
$\nu(x,z,t)\,\leq\,\nu(x,y,t)\,\diamond\,\nu(y,z,t).$\\
Then $(X,\mu,\nu,\ast,\diamond)$ is called non-Archimedean intuitionistic fuzzy metric space.
\end{remark}

\begin{definition}\cite{abdul}
Let $(X,\mu,\nu,\ast,\diamond)$ be an intuitionistic fuzzy metric space. A mapping $f:X\rightarrow\,X$ is intuitionistic fuzzy contractive if there exists $k\in (0,1)$ such that $\frac{1}{\mu\left(f(x),f(y),t\right)}-1\,\leq\,k\,\left(\frac{1}{\mu(x,y,t)}-1\right)\;$ and $\;\frac{1}{\nu\left(f(x),f(y),t\right)}-1\,\leq\,\frac{1}{k}\,\left(\frac{1}{\nu(x,y,t)}-1\right)\;$ for all $x,y\in X$ and $t>0.\;$ ($k$ is called contractive constant of $f$.)
\end{definition}

\begin{definition}\cite{abdul}
Let $(X,\mu,\nu,\ast,\diamond)$ be an intuitionistic fuzzy metric space. We will say that the sequence $\{x_n\}_n$ in $X$ is intuitionistic fuzzy contractive if there exists $k\in (0,1)$ such that $\frac{1}{\mu\left(x_{n+1},x_{n+2},t\right)}-1\,\leq\,k\,\left(\frac{1}{\mu(x_n,x_{n+1},t)}-1\right)\;$ and $\;\frac{1}{\nu\left(x_{n+1},x_{n+2},t\right)}-1\,\leq\,\frac{1}{k}\,\left(\frac{1}{\nu(x_{n},x_{n+1},t)}-1\right)\;$ for all $t>0$ and $n\in\mathbb{N}.$
\end{definition}
\[\]

%=======================================================================
\section{\bf Intuitionistic fuzzy $\Psi$-$\Phi$-contractive mappings}
%=======================================================================
\begin{definition}
Let $(X,\mu,\nu,\ast,\diamond)$ be an intuitionistic fuzzy metric space.\\\\
(i) A sequence $\{x_n\}_n$ in $X$ is called M-Cauchy sequence, if for each $\epsilon\in (0,1)$ and $t>0$ there exists $n_0\in \mathbb{N}$ such that $\mu(x_n,x_m,t)>1-\epsilon$ and $\nu(x_n,x_m,t)<\epsilon\,$ for all $m,n\geq\,n_0.\\\\$
(ii) A sequence $\{x_n\}_n$ in $X$ is called G-Cauchy sequence if $\,\mathop {\lim }\limits_{t\;\, \to
\,\;\infty } \,\,\mu\,\left(\,x_n,x_{n+m},t\, \right)=1\,$ and $\,\mathop {\lim }\limits_{t\;\, \to
\,\;\infty} \,\,\nu\,\left(\,x_n,x_{n+m},t\, \right)=0 $ for each $m\in \mathbb{N}$ and $t>0$.
\end{definition}

\begin{definition}
A sequence $\{x_n\}_n$ in an intuitionistic fuzzy metric space $(X,\mu,\nu,\ast,\diamond)$ is said to converge to $x\in X \;$ if $\;\mathop {\lim }\limits_{t\;\, \to
\,\;\infty } \,\,\mu\,\left(\,x_n,x,t\, \right)=1\,$ and $\,\mathop {\lim }\limits_{t\;\, \to
\,\;\infty} \,\,\nu\,\left(\,x_n,x,t\, \right)=0 $ for all $t>0$.
\end{definition}

\begin{definition}
Let $\Psi$ be the class of all mappings $\psi:[0,1]\rightarrow\,[0,1]$ such that $\psi$ is continuous, non-increasing and $\psi(t)<t,\,\forall\,t\in(0,1).\;\,$ Let $\Phi$ be the class of all mappings $\phi:[0,1]\rightarrow\,[0,1]$ such that $\phi$ is continuous, non-decreasing and $\phi(t)>t,\,\forall\,t\in(0,1).\;\,$ Let $(X,\mu,\nu,\ast,\diamond)$ be an intuitionistic fuzzy metric space and $\psi\in \Psi$ and $\phi\in\Phi.\,$ A mapping $f:X\rightarrow\,X$ is called an intuitionistic fuzzy $\psi$-$\phi$-contractive mapping if the following implications hold:
\[\mu(x,y,t)>0\,\Rightarrow\;\psi\left(\mu\left(f(x),f(y),t\right)\right)\,\geq\,\mu(x,y,t)\]
\[\nu(x,y,t)<1\,\Rightarrow\;\phi\left(\nu\left(f(x),f(y),t\right)\right)\,\leq\,\nu(x,y,t).\]
\end{definition}

\begin{example}
Let $(X,\mu,\nu,\ast,\diamond)$ be an intuitionistic fuzzy metric space and $f:X\rightarrow\,X$ satisfies
$\frac{1}{\mu\left(f(x),f(y),t\right)}-1\,\leq\,k\,\left(\frac{1}{\mu(x,y,t)}-1\right)\;$ and $\;\frac{1}{\nu\left(f(x),f(y),t\right)}-1\,\leq\,\frac{1}{k}\,\left(\frac{1}{\nu(x,y,t)}-1\right)\;$ for all $x,y\in X$ and $t>0.\;$ Then for each $k\in\,(0,1)$, $f$ is an intuitionistic fuzzy $\psi$-$\phi$-contractive mapping, with
\[\psi(t)=\,-k,\;\;\;\phi(t)=\frac{t}{(1-k)t+k}\]
\end{example}

\begin{example}
Let $X$ be a non-empty set with at least two elements. If we define the fuzzy set $(X,\,\mu,\,\nu)$ by $\mu(x,x,t)=1$ and $\nu(x,x,t)=0$ for all $x\in X$ and $t>0;\,$ and \\
$\mu\,(\,x\,,\,t\,)\;=\begin{cases}0,\;\;\;\,\hspace{0.5 cm}if\;t\,\leq\,1\hspace{0 cm}
\\ 1,\;\;\;\hspace{0.5 cm}if\,\;t\,>\,1\hspace{-2.5 cm}\end{cases}\hspace{2 cm}$
$\nu\,(\,x\,,\,t\,)\;=\begin{cases}1,\;\;\;\,\hspace{0.5 cm}if\;t\,\leq\,1\hspace{0 cm}
\\ 0,\;\;\;\hspace{0.5 cm}if\,\;t\,>\,1\hspace{-2.5 cm}\end{cases}$\\\\
for all $x,y\in X,\;x\neq\,y,$ then $(X,\mu,\nu,\ast,\diamond)$ is an M-complete non-Archimedean fuzzy metric space under any continuous t-norm $\ast$ and continuous t-conorm $\diamond$. Now,\\
$\mu(x,y,t)>0\;\\\Rightarrow\,\mu(x,y,t)=1\;\\\Rightarrow\,\psi\left(\mu(f(x),f(y),t)\right)\geq\,\mu(x,y,t)=1 \;\\\Rightarrow\,\psi\left(\mu(f(x),f(y),t)\right)=1=\,\mu(x,y,t);\;$ and\\
$\nu(x,y,t)<1\;\\\Rightarrow\,\nu(x,y,t)=0\;\\\Rightarrow\,\phi\left(\nu(f(x),f(y),t)\right)\leq\,\nu(x,y,t)=0 \;\\\Rightarrow\,\phi\left(\nu(f(x),f(y),t)\right)=0=\,\nu(x,y,t).\;$\\
Therefore every mapping $f:X\rightarrow\,X$ is an intuitionistic fuzzy $\psi$-$\phi$-contractive mapping.
\end{example}

\begin{definition}
An intuitionistic fuzzy $\psi$-$\phi$-contractive sequence in an intuitionistic fuzzy metric space $(X,\mu,\nu,\ast,\diamond)$ is any sequence $\{x_n\}_n$ in $X$ such that
\[\psi\left(\mu\left(x_{n+1},x_{n+2},t\right)\right)\,\geq\,\mu(x_{n+1},x_n,t)\]
\[\phi\left(\nu\left(x_{n+1},x_{n+2},t\right)\right)\,\leq\,\nu(x_{n+1},x_n,t).\]
An intuitionistic fuzzy metric space $(X,\mu,\nu,\ast,\diamond)$ is called M-complete (G-complete) if every M-Cauchy (G-Cauchy) sequence is convergent in $X$.
\end{definition}

%=====================================
\section{\bf Fixed point theorems}
%=====================================
\begin{theorem}
Let $(X,\mu,\nu,\ast,\diamond)$ be an M-complete non-Archimedean intuitionistic fuzzy metric space and $f:X\rightarrow\,X$ be an intuitionistic fuzzy $\psi$-$\phi$-contractive mapping. If there exists $x\in X$ such that $\mu\left(x,f(x),t\right)>0$ and $\nu\left(x,f(x),t\right)<1$ for all $t>0,\;$ then $f$ has a unique fixed point.
\end{theorem}
{\bf Proof.} Let $x\in X$ be such that $\mu(x,f(x),t)>0$ and $\nu(x,f(x),t)<1,\;t>0\,$ and $x_n=f^n(x),\;n\in\mathbb{N}$, we have for all $t>0\\$
\[\mu(x_0,x_1,t)\leq\,\psi(\mu(x_1,x_2,t))<\mu(x_1,x_2,t)\]
\[\nu(x_0,x_1,t)\geq\,\phi(\nu(x_1,x_2,t))>\nu(x_1,x_2,t)\]
and
\[\mu(x_1,x_2,t)\leq\,\psi(\mu(x_2,x_3,t))<\mu(x_2,x_3,t)\]
\[\nu(x_1,x_2,t)\geq\,\phi(\nu(x_2,x_3,t))>\nu(x_2,x_3,t)\]
Hence by induction $\forall\,t>0,\;\mu(x_n,x_{n+1},t)<\mu(x_{n+1},x_{n+2},t)$ and $\nu(x_n,x_{n+1},t)>\nu(x_{n+1},x_{n+2},t).\;$ Therefore, for every $t>0,\;\{\mu(x_n,x_{n+1},t)\}$ is a non-increasing sequence of numbers in $(0,1]$ and $\{\nu(x_n,x_{n+1},t)\}$ is a non-decreasing sequence of numbers in $[0,1).\\$
Fix $t>0$. Denote $\mathop {\lim }\limits_{n\;\to\;\infty}\,\mu(x_n, x_{n+1},t) $ by $l$ and $\mathop {\lim }\limits_{n\;\to\;\infty}\,\nu(x_n, x_{n+1},t) $ by $m$. Then we have $l\in\,[\,0,1\,]$ and $m\in\,[\,0,1\,]$. Since  $\psi\left(\mu\left(x_{n+1},x_{n+2},t\right)\right)\,\geq\,\mu(x_{n},x_{n+1},t)$ and $\psi$ is continuous, $\psi(l)\geq\,l$. This implies $l=1$. Also, since $\phi\left(\nu\left(x_{n+1},x_{n+2},t\right)\right)\,\leq\,\nu(x_{n+1},x_n,t)$ and $\phi$ is continuous, $\phi(m)\leq\,m.$ This implies $m=0$. Therefore,\\
$\mathop {\lim }\limits_{n\;\to\;\infty}\,\mu(x_n, x_{n+1},t)\,=\,1$ and $\mathop {\lim }\limits_{n\;\to\;\infty}\,\nu(x_n, x_{n+1},t)\,=\,0.\\$
If $\{x_n\}_n$ is not a M-cauchy sequence then there are $\epsilon\,\in\,(0,1)$ and $t>0$ such that for each $k\in\,\mathbb{N}$ there exist $m(k),\,n(k)\in\,\mathbb{N}\;$ with $\;m(k)>n(k)\geq\,k$ and \\
$\mu(x_{m(k)},\,x_{n(k),t})\,\leq\,1-\epsilon\;$ and $\;\nu(x_{m(k)},\,x_{n(k),t})\,\geq\,\epsilon\\$
Let for each $k,\;m(k)$ be the least positive integer exceeding $n(k)$ satisfying the above property, that is,\\
$\mu(x_{m(k)-1},\,x_{n(k)})\,\geq\,1-\epsilon$ and $\mu(x_{m(k)},\,x_{n(k)})\,\leq\,1-\epsilon\,.$ Also,\\
$\nu(x_{m(k)-1},\,x_{n(k)})\,\leq\,\epsilon$ and $\nu(x_{m(k)},\,x_{n(k)})\,\geq\,\epsilon\,.\\$
Then for each positive integer $k$,
\[1-\epsilon\,\geq\,\mu(x_{m(k)},\,x_{n(k)},t)\hspace{8 cm}\]
\[\geq\,\,\mu(x_{m(k)-1},\,x_{n(k)},t)\,\ast\,\mu(x_{m(k)-1},\,x_{m(k)},t)\hspace{2.5 cm}\]
\[\geq\,\,(1-\epsilon)\,\ast\,\mu(x_{m(k)-1},\,x_{m(k)},t).\hspace{4.5 cm}\]
and \[\epsilon\,\leq\,\nu(x_{m(k)},\,x_{n(k)},t)\hspace{7 cm}\]
\[\leq\,\,\nu(x_{m(k)-1},\,x_{n(k)},t)\,\diamond\,\nu(x_{m(k)-1},\,x_{m(k)},t)\hspace{2.25 cm}\]
\[\leq\,\,\epsilon\,\diamond\,\nu(x_{m(k)-1},\,x_{m(k)},t).\hspace{5 cm}\]
Taking limit as $k\rightarrow\,\infty$ we have,\\
$\;\mathop {\lim }\limits_{k\;\to\;\infty}\,\{(1-\epsilon)\,\ast\,\mu(x_{m(k)-1}\,,\, x_{m(k)},t)\}\,=\,(1-\epsilon)\,\ast\,\mathop {\lim }\limits_{k\;\to\;\infty}\,\mu(x_{m(k)-1}\,,\, x_{m(k)},t)\,=\,(1-\epsilon)\,\ast\,1\,=\,(1-\epsilon)\;$ and \\$\;\mathop {\lim }\limits_{k\;\to\;\infty}\,\{\epsilon\,\diamond\,\nu(x_{m(k)-1}\,,\, x_{m(k)},t)\}\,=\,\epsilon\,\diamond\,\mathop {\lim }\limits_{k\;\to\;\infty}\,\nu(x_{m(k)-1}\,,\, x_{m(k)},t)\,=\,\epsilon\,\diamond\,0\,=\,\epsilon\,$. \\It follows that
$\;\,\mathop {\lim }\limits_{k\;\to\;\infty}\,\mu(x_{m(k)},\,x_{n(k)},t)\,=\,1-\epsilon$ and $\mathop {\lim }\limits_{k\;\to\;\infty}\,\mu(x_{m(k)},\,x_{n(k)},t)\,=\,\epsilon$.\\
Now, $\mu(x_{m(k)},\,x_{n(k)},t)\,\leq\,\psi\left(\mu(x_{m(k)+1},\,x_{n(k)+1},t)\right)$ and $\nu(x_{m(k)},\,x_{n(k)},t)\,\geq\,\phi\left(\nu(x_{m(k)+1},\,x_{n(k)+1},t)\right)\,.\;$ Since $\psi$  and $\phi$ are continuous taking limit as $k\rightarrow\,\infty$ we have,\\
$1-\epsilon\,\leq\,\psi(1-\epsilon)\,<\,1-\epsilon$ and $\epsilon\,\geq\,\phi(\epsilon)\,>\,\epsilon\,,\;$ which are contradictions. Thus $\{x_n\}_n$  is a M-cauchy sequence.\\
If $\,\mathop {\lim }\limits_{k\;\to\;\infty}\,x_n\,=\,y\;$ then from $\psi\left(\mu(f(y),f(x_n),t)\right)\,\geq\,\mu(y,x_n,t)$ and $\phi\left(\nu(f(y),f(x_n),t)\right)\,\leq\,\nu(y,x_n,t)\;$ it follows that $\;x_{n+1}\rightarrow\,f(y).\\$
Therefore we have \\
$\mu(y,f(y),t)\,\geq\,\mu(y,x_n,t)\,\ast\,\mu(x_n,x_{n+1},t)\,\ast\,\mu(x_{n+1},f(y),t)\,\rightarrow\,1\;$ as $n\rightarrow\,\infty.$
This implies $\mu(y,f(y),t)\,=\,1.\\$
$\nu(y,f(y),t)\,\leq\,\nu(y,x_n,t)\,\diamond\,\nu(x_n,x_{n+1},t)\,\diamond\,\nu(x_{n+1},f(y),t)\,\rightarrow\,0\;$ as $n\rightarrow\,\infty.$
This implies $\nu(y,f(y),t)\,=\,0.$
Hence, $f(y)\,=\,y.$\\
If $x,y$ are fixed points of $f$ then $\mu(f(x),f(y),t)=\mu(x,y,t)\leq\,\psi(\mu(f(x),f(y),t))$ and $\nu(f(x),f(y),t)=\nu(x,y,t)\geq\,\phi(\nu(f(x),f(y),t)),\;\forall\,t>0.\\$
If $x\neq y\;$ then $\mu(x,y,s)<1$ and $\nu(x,y,s)>0$ for some $s>0\;$ i.e., $\;0<\mu(x,y,s)<1$ and $0<\nu(x,y,s)<1$ hold, impllying\\
$\mu\left(f(x),f(y),s\right)\leq\psi(\mu\left(f(x),f(y),s\right))<\mu(f(x),f(y),s)$ and $\nu\left(f(x),f(y),s\right)\geq\phi(\nu\left(f(x),f(y),s\right))>\nu(f(x),f(y),s),\;$ which are contradictions.\\
Thus $\;x=y\,.$\\
This completes the proof.

\begin{lemma}\label{l1}
Let $(X,\mu,\nu,\ast,\diamond)$ be a non-Archimedean intuitionistic fuzzy metric space. If $\{x_n\}_n$ and $\{y_n\}_n$ be two sequences in $X$ converges to $x$ and $y$ respectively then $\mathop {\lim }\limits_{n\;\, \to
\,\;\infty} \,\,\mu\,\left(\,x_n,y_n,t\, \right)=\mu\,(x,y,t)\,$ and $\mathop {\lim }\limits_{n\;\, \to
\,\;\infty} \,\,\nu\,\left(\,x_n,y_n,t\, \right)=\nu\,(x,y,t).\,$
\end{lemma}
{\bf Proof.} Since $(X,\mu,\nu,\ast,\diamond)$ be a non-Archimedean intuitionistic fuzzy metric space, therefore\\
$\mu\,\left(\,x_n,y_n,t\, \right)\,\geq\,\mu(x_n,x,t)\,\ast\,\mu(x,y,t)\,\ast\,\mu(y,y_n,t)\\$
$\Rightarrow\;\mathop {\lim }\limits_{n\;\, \to\;\infty}\mu\,\left(\,x_n,y_n,t\, \right)\,\geq\,1\,\ast\,\mu(x,y,t)\,\ast\,1\;=\,\mu(x,y,t)\,.$\\
and $\nu\,\left(\,x_n,y_n,t\, \right)\,\leq\,\nu(x_n,x,t)\,\diamond\,\nu(x,y,t)\,\diamond\,\nu(y,y_n,t)\\$
$\Rightarrow\;\mathop {\lim }\limits_{n\;\, \to\;\infty}\nu\,\left(\,x_n,y_n,t\, \right)\,\leq\,0\,\diamond\,\nu(x,y,t)\,\diamond\,0\;=\,\nu(x,y,t)\,.$\\\\
Also, $\mu(x,y,t)\,\geq\,\mu(x,x_n,t)\,\ast\,\mu(x_n,y_n,t)\,\ast\,\mu(y,y_n,t)\\$
$\Rightarrow\;\mu(x,y,t)\,\geq\,1\,\ast\,\mathop {\lim }\limits_{n\;\, \to
\,\;\infty} \mu(x_n,y_n,t)\,\ast\,1\,=\,\mathop {\lim }\limits_{n\;\, \to
\,\;\infty} \mu(x_n,y_n,t)\\$
and $\nu(x,y,t)\,\geq\,\nu(x,x_n,t)\,\ast\,\nu(x_n,y_n,t)\,\ast\,\nu(y,y_n,t)\\$
$\Rightarrow\,\nu(x,y,t)\,\geq\,0\,\diamond\,\mathop {\lim }\limits_{n\;\, \to
\,\;\infty} \nu(x_n,y_n,t)\,\diamond \,0\,=\,\mathop {\lim }\limits_{n\;\, \to
\,\;\infty} \nu(x_n,y_n,t).\\$
Hence the proof.

\begin{theorem}
Let $(X,\mu,\nu,\ast,\diamond)$ be a compact non-Archimedean intuitionistic fuzzy metric space. Let $f:X\rightarrow\,X$ be an intuitionistic fuzzy $\psi$-$\phi$-contractive mapping. Then $f$ has a unique fixed point.
\end{theorem}
{\bf Proof.} Let $x\in X$ and $x_n=f^n (x),\;n\in\mathbb{N}.\;$ Assume $x_n\neq x_{n+1}$ for each $n$ (if not $f(x_n)=x_n$). \\
Now assume $x_n\neq x_m\;(n\neq m)$, otherwise for $m<n$ we get\\
$\mu(x_n,x_{n+1},t)=\mu(x_m,x_{m+1},t)\leq \psi\left(\mu(x_{m+1},x_{m+2},t)\right)<\mu(x_{m+1},x_{m+2},t)<\cdots<\mu(x_n,x_{n+1},t)\;$ and \\
$\nu(x_n,x_{n+1},t)=\nu(x_m,x_{m+1},t)\geq \phi\left(\nu(x_{m+1},x_{m+2},t)\right)>\nu(x_{m+1},x_{m+2},t)>\cdots>\nu(x_n,x_{n+1},t),\;$ a contradiction.\\
Since $X$ is compact, $\{x_n\}_n$ in $X$ has a convergent subsequence $\{x_{n_i}\}_{i\in\mathbb{N}}$ (say). Let $\{x_{n_i}\}_{i\in\mathbb{N}}$ converges to $y$. We also assume that $y,f(y)\notin \{x_n:n\in \mathbb{N}\}$ (if not, choose a subsequence with such a property). According to the above assumptions we may now write for all $i\in\mathbb{N}$ and $t>0$
\[\mu(x_{n_i},y,t)\leq \psi(\mu\left(f(x_{n_i}),f(y),t\right))<\mu\left(f(x_{n_i}),f(y),t\right)\]
\[\nu(x_{n_i},y,t)\geq \phi(\nu\left(f(x_{n_i}),f(y),t\right))>\nu\left(f(x_{n_i}),f(y),t\right)\]
Since $\psi$ and $\phi$ are continuous for all $x,y\in X$. From lemma 4.2 we obtain\\
$\mathop {\lim }\limits_{i\;\, \to\;\infty}\,\mu(x_{n_i},y,t)\leq\,\mathop {\lim }\limits_{i\;\, \to\;\infty}\,\mu\left(f(x_{n_i}),f(y),t\right)\\$
$\Rightarrow\;1\leq\,\mathop {\lim }\limits_{i\;\, \to\;\infty}\,\mu\left(f(x_{n_i}),f(y),t\right)\\$
$\Rightarrow\;\mathop {\lim }\limits_{i\;\, \to\;\infty}\,\mu\left(f(x_{n_i}),f(y),t\right)=1.$ and\\
$\mathop {\lim }\limits_{i\;\, \to\;\infty}\,\nu(x_{n_i},y,t)\geq\,\mathop {\lim }\limits_{i\;\, \to\;\infty}\,\nu\left(f(x_{n_i}),f(y),t\right)\\$
$\Rightarrow\;0\geq\,\mathop {\lim }\limits_{i\;\, \to\;\infty}\,\nu\left(f(x_{n_i}),f(y),t\right)\\$
$\Rightarrow\;\mathop {\lim }\limits_{i\;\, \to\;\infty}\,\nu\left(f(x_{n_i}),f(y),t\right)=0.$ i.e.,
\begin{equation}\label{e1}
f(x_{n_i})\; \rightarrow\; f(y)
\end{equation}
Similarly, we obtain\\
\begin{equation}\label{e2}
f^2(x_{n_i})\; \rightarrow\; f^2(y)
\end{equation}
Now, we see that\\
$\mu(x_{n_1},f(x_{n_1}),t)\leq \psi(\mu\left(f(x_{n_1}),f^2(x_{n_1}),t\right))< \mu\left(f(x_{n_1}),f^2(x_{n_1}),t\right)<\cdots<\mu(x_{n_i},f(x_{n_i}),t) < \mu\left(f(x_{n_i}),f^2(x_{n_i}),t\right)<\cdots<1.$ and\\
$\nu(x_{n_1},f(x_{n_1}),t)\geq \phi(\nu\left(f(x_{n_1}),f^2(x_{n_1}),t\right))> \nu\left(f(x_{n_1}),f^2(x_{n_1}),t\right)>\cdots>\nu(x_{n_i},f(x_{n_i}),t) > \nu\left(f(x_{n_i}),f^2(x_{n_i}),t\right)>\cdots>0.$\\
Thus $\{\mu(x_{n_i},f(x_{n_i}),t)\}_{i\in\mathbb{N}}$ and $\{\mu\left(f(x_{n_i}),f^2(x_{n_i}),t\right)\}_{i\in\mathbb{N}}$ converges to a common limit. Also, $\{\nu(x_{n_i},f(x_{n_i}),t)\}_{i\in\mathbb{N}}$ and $\{\nu\left(f(x_{n_i}),f^2(x_{n_i}),t\right)\}_{i\in\mathbb{N}}$ converges to a common limit.\\
So, by (\ref{e1}), (\ref{e2}) and lemma \ref{l1} we get\\
$\mu(y,f(y),t)=\mu\left(\mathop {\lim }\limits_{i\;\, \to\;\infty}\,x_{n_i}, \,f(\mathop {\lim }\limits_{i\;\, \to\;\infty} \,x_{n_i}),\,t\right) =
\mathop {\lim }\limits_{i\;\, \to\;\infty}\, \mu\left(\,x_{n_i}, \,f(x_{n_i}),\,t\right) \\= \mathop {\lim }\limits_{i\;\, \to\;\infty}\, \mu\left(\,f(x_{n_i}), \,f^2(x_{n_i}),\,t\right)=\mu\left(\,f(\mathop {\lim }\limits_{i\;\, \to\;\infty}\,x_{n_i}), \,f^2(\mathop {\lim }\limits_{i\;\, \to\;\infty}\,x_{n_i}),\,t\right)=\mu(f(y),f^2(y),t)$ and \\
$\nu(y,f(y),t)=\nu\left(\mathop {\lim }\limits_{i\;\, \to\;\infty}\,x_{n_i}, \,f(\mathop {\lim }\limits_{i\;\, \to\;\infty} \,x_{n_i}),\,t\right) =
\mathop {\lim }\limits_{i\;\, \to\;\infty}\, \nu\left(\,x_{n_i}, \,f(x_{n_i}),\,t\right) \\= \mathop {\lim }\limits_{i\;\, \to\;\infty}\, \nu\left(\,f(x_{n_i}), \,f^2(x_{n_i}),\,t\right)=\nu\left(\,f(\mathop {\lim }\limits_{i\;\, \to\;\infty}\,x_{n_i}), \,f^2(\mathop {\lim }\limits_{i\;\, \to\;\infty}\,x_{n_i}),\,t\right)=\nu(f(y),f^2(y),t)$ for all $t>0.$\\
Suppose $f(y)\neq y$, then we have $\mu(y,f(y),t)\leq \psi(\mu\left(f(y),f^2(y),t\right))<\mu\left(f(y),f^2(y),t\right)$ and $\nu(y,f(y),t)\geq \phi(\nu\left(f(y),f^2(y),t\right))>\nu\left(f(y),f^2(y),t\right),\;t>0,\;$ a contradiction.\\
Hence $y=f(y)$ is a fixed point.\\
If $x,y$ are fixed points of $f$ then $\mu(f(x),f(y),t)=\mu(x,y,t)\leq\,\psi(\mu(f(x),f(y),t))$ and $\nu(f(x),f(y),t)=\nu(x,y,t)\geq\,\phi(\nu(f(x),f(y),t)),\;\forall\,t>0.\\$
Suppose that $x\neq y,\;$ then $\mu(x,y,s)<1$ and $\nu(x,y,s)>0$ for some $s>0$ i.e., $0<\mu(x,y,s)<1$ and $0<\nu(x,y,s)<1$ hold, impllying\\
$\mu\left(f(x),f(y),s\right)\leq\psi(\mu\left(f(x),f(y),s\right))<\mu(f(x),f(y),s)$ and $\nu\left(f(x),f(y),s\right)\geq\phi(\nu\left(f(x),f(y),s\right))>\nu(f(x),f(y),s),\;$ which are contradictions.\\
Therefore it must be the case that $x=y.\\$
Hence the proof.

\end{document}